\documentclass[10pt]{amsart}
\usepackage{tikz,array, verbatim}
\usepackage{amsfonts, amsmath, latexsym, epsfig, caption, cite}
\usepackage{amssymb, color}

\usepackage{epsf}
\usepackage{array}
\usepackage{ragged2e}
\usepackage{hyperref}
\usepackage{longtable}
\usepackage{pdflscape}
\usepackage{booktabs}
\newcolumntype{C}[1]{>{\centering\let\newline\\\arraybackslash\hspace{0pt}}m{#1}}
\newcolumntype{L}[1]{>{\raggedright\let\newline\\\arraybackslash\hspace{0pt}}m{#1}}
\newcolumntype{R}[1]{>{\raggedleft\let\newline\\\arraybackslash\hspace{0pt}}m{#1}}

\newif\iflandscapetable
\landscapetabletrue

\newif\ifshowvc
\showvctrue
\showvcfalse  

\ifshowvc
\input{vc}
\fi

\DeclareMathOperator{\SL}{SL}

\DeclareMathOperator{\GL}{GL}
\DeclareMathOperator{\PSL}{PSL}

\DeclareMathOperator{\Eis}{Eis}
\DeclareMathOperator{\cusp}{cusp}

\DeclareMathOperator{\Vor}{Vor}

\DeclareMathOperator{\Cl}{Cl}
\DeclareMathOperator{\logtor}{logtor}

\newcommand{\cS}{\ensuremath{\mathcal{S}}}

\newcommand{\OO}{\ensuremath{\mathcal{O}}}

\newcommand{\QQ}{\ensuremath{\mathbb{Q}}}
\newcommand{\CC}{\ensuremath{\mathbb{C}}}
\newcommand{\ZZ}{\ensuremath{\mathbb{Z}}}

\DeclareMathOperator{\tors}{\mathrm{tors}}

\theoremstyle{plain}
\newtheorem{theorem}{Theorem}[section]

\newtheorem{results}[theorem]{Results}
\theoremstyle{definition}

\theoremstyle{remark}

\def\QuotS#1#2{\leavevmode\kern-.0em\raise.2ex\hbox{$#1$}\kern-.1em/\kern-.1em\lower.25ex\hbox{$#2$}}

\newcommand{\fH}{\mathfrak{H}}

\begin{document}
\author[H. Gangl]{Herbert Gangl}
\address{H. Gangl, Department of Mathematical Sciences, Stockton Road, Durham DH1 3LE, United Kingdom}
\email{herbert.gangl@durham.ac.uk}
\urladdr{\url{http://maths.dur.ac.uk/~dma0hg/}}

\author[P. E. Gunnells]{Paul E. Gunnells}
\address{P. E. Gunnells, Department of Mathematics and Statistics, LGRT 1115L, University of Massachusetts, Amherst, MA 01003, USA}
\email{gunnells@math.umass.edu}
\urladdr{\url{https://www.math.umass.edu/~gunnells/}}

\author[J. Hanke]{Jonathan Hanke}
\address{J. Hanke, Princeton, NJ 08542, USA}
\email{jonhanke@gmail.com}
\urladdr{\url{http://www.jonhanke.com}}
\author[D. Yasaki]{Dan Yasaki}
\address{D. Yasaki, Department of Mathematics and Statistics,
University of North Carolina at Greensboro, Greensboro, NC 27412, USA}
\email{d\_yasaki@uncg.edu}
\urladdr{\url{http://www.uncg.edu/~d_yasaki/}}

\thanks{PG was partially supported by the NSF under contracts DMS
1101640 and DMS 1501832.  DY was partially supported by the Simons
Foundation under contract 848154.  The authors thank the American
Institute of Mathematics and the MPIM Bonn for their hospitality and
excellent working conditions. We also thank G.~Harder for inspiring
discussions and J.~Schwermer for helpful comments.}



\title{On the cohomology of $\GL_{2}$ and $\SL_{2}$ over imaginary
quadratic fields} 

\date{\today}

\begin{abstract}
We report on computations of the cohomology of $\GL_{2} (\OO_{D})$ and
$\SL_{2} (\OO_{D})$, where $D<0$ is a fundamental discriminant.  These
computations go well beyond earlier results of Vogtmann
\cite{vogtmann} and Scheutzow \cite{scheutzow}.  We use the technique
of homology of Voronoi complexes, and our computations recover the
integral cohomology away from the primes $2, 3$.  We observed
exponential growth in the torsion subgroup of $H^{2}$ as $|D|$
increases, and compared our data to bounds of Rohlfs \cite{rohlfs.bound}.
\end{abstract}

\maketitle
\ifshowvc
\let\thefootnote\relax
\footnotetext{Base revision~\GITAbrHash, \GITAuthorDate,
\GITAuthorName.}
\fi

\section{Introduction}

\subsection{} Let $F$ be an imaginary quadratic field of discriminant
$D<0$, let $\OO = \OO_{D}$ be its ring of integers, and let $\Gamma$
be the group $\GL_{2} (\OO)$ or $\SL_{2} (\OO)$.  The purpose of this
paper is to computationally investigate the group cohomology of
$\Gamma$.  Such computations for $\Gamma$ and its congruence subgroups
have been investigated by many authors,
cf.~\cite{rahm1,rahm2,berkove,schwermer-vogtmann,sengun,cremona,yasaki,mendoza};
of particular relevance to our results is the prior work of Vogtmann
\cite{vogtmann} and Scheutzow \cite{scheutzow}.

In this paper we significantly extend this prior work by using
Voronoi's explicit reduction theory \cite{aim-imagquad} to compute the
integral cohomology of $\Gamma $ modulo small primes and for a large
range of discriminants. More precisely, the only primes dividing the
orders of the finite subgroups of $\Gamma$ for the groups under
consideration are $2$ and $3$.
There is a (homologically-indexed) complex $\Vor_{D} = (V_{*}
(\Gamma), d_{*})$ of $\ZZ [\Gamma]$-modules, called
the \emph{Voronoi complex}, such that 
\[
H_{i} (V_{*} (\Gamma)) \simeq_{{3}} H^{3-i} (\Gamma ; \ZZ), \quad i=1,2,3.
\]
Here the notation $\simeq_{n}$ means the following.  For any positive
integer $n$, let $\cS_{n}$ be the Serre class of all finite abelian
groups with orders only divisible by primes $\leq n$.  Then for two
finitely generated abelian groups $A, B$ the notation $A \simeq_{n} B$
means that $A$ is isomorphic to $B$ modulo $\cS_{n}$.  For more
details about the construction of the complex $\Vor_{D}$ and the
connection with group cohomology, we refer to \cite{aim-imagquad}.

Our main experimental results and observations are as follows:

\begin{results}
\leavevmode
\begin{enumerate}
\item We computed the Voronoi complex and its homology for
$\GL_2(\OO_{D})$ with $-D\leq 2099$ and for $\SL_2(\OO_{D})$ with
$-D\leq 1247$, and thus computed the cohomology with integral
coefficients (away from the primes $2, 3$) for these groups.
\item We understand the Betti numbers explicitly.  In particular we
see how the Eisenstein cohomology contributes.  For $\SL_{2} $ the
pattern of Eisenstein cohomology agrees with the description in
\cite{vogtmann} for the Eisenstein classes in the group homology (\S \ref{ss:ccec}).
\item Away from the primes $2$ and $3$, we observe that the order of
the torsion subgroup in the Voronoi homology $H_{1}$ is a square.
Away from $2$, it is typically the same for $\GL_{2}$ and $\SL_{2}$.
There are however some exceptions, such as $D=-1151$ (since $31$
divides the order of the torsion subgroup for $\SL_{2}$ but not
$\GL_{2}$).  Apart from such examples which account for approximately
3\% of our data, the main difference lies in the order of the
$2$-torsion.
\item For $\GL_{2}$, we observe exponential growth in the order of the
torsion subgroup of the Voronoi homology $H_{1}$ as $|D|$ increases (\S\ref{ss:eg});
Exponential growth for $\SL_{2}$ is not apparent from our data.
This includes not only growth in the order of the group but also the
appearance of large prime factors (such as $14116228597231$ for
$\GL_{2} (\OO_{-2087})$).  We also noted that the torsion in the
Voronoi homology group $H_{2}$ for either $\GL_{2}/\SL_{2}$ is killed
by $12$.
\item We compare our data to a lower bound on the dimension of the
cuspidal cohomology in $H^1$ for $\PSL_{2} (\OO)$ due to Rohlfs
\cite{rohlfs.bound} and observe in many cases the bound is sharp
(Figure~\ref{fig:rohlfs} in \S\ref{ss:rbcf}).
\item We observe surprising regularity in the linear growth of the minimal number of
(infinite and finite) cyclic factors in the Voronoi homology $H_1$ for $\GL_2(\OO_{D})$ as
$|D|$ increases (Figure~\ref{fig:Z} in \S\ref{ss:rbcf}).
\end{enumerate}
\end{results}

\section{Background and structure of the cohomology}

\subsection{Complex group cohomology and automorphic forms} First we
consider the Voronoi homology with $\CC$-coefficients $H_{*}(V_{*}
(\Gamma))\otimes \CC$.  We have
\[
H_{i}(V_{*} (\Gamma))\otimes \CC \simeq H^{3-i} (\Gamma ; \CC).
\]
It is known that the complex group cohomology $H^{*} (\Gamma ; \CC)$
can be computed in terms of certain automorphic forms for $\Gamma$;
for a general discussion of this interpretation we refer to
\cite{vogan,li.schwermer.survey}.

Furthermore the symmetric space $\fH_{3}$ for $\Gamma$ has a partial
compactification $\overline{\fH}_{3}$ due to Borel--Serre
\cite{borel.serre}, which leads to a structural decomposition of the
cohomology.  The quotient $\Gamma \backslash \overline{\fH}_{3}$ is a
compact three-dimensional orbifold, and we have
\[
H^{*} (\Gamma ;
\CC) \simeq H^{*} (\Gamma \backslash \fH_{3}; \CC )
\simeq H^{*} (\Gamma \backslash \overline{\fH}_{3}; \CC ).
\]
The inclusion of the boundary $\partial (\Gamma \backslash
\overline{\fH}_{3})$ in $\Gamma \backslash \overline{\fH}_{3}$ induces
a restriction map
\[
r\colon H^{*} (\Gamma \backslash
\overline{\fH}_{3}) \rightarrow H^{*} (\partial (\Gamma \backslash
\overline{\fH}_{3})).
\]
The kernel $H^{*}_{!}$ of $r$ is called the \emph{interior
cohomology}.  There is a complement $H^{*}_{\Eis}$ of $H^{*}_{!}$ in
$H^{*} (\Gamma ;\CC)$ called the \emph{Eisenstein cohomology}; its
construction, which uses Eisenstein series, is due to Harder
\cite{harder.gl2,harder-schwermerconf,harder.budapest}.  For the
$\Gamma$ under consideration, the interior cohomology can be further 
identified with the \emph{cuspidal cohomology} $H^{*}_{\cusp}$ \cite{grunewald.schwermer.2}; by
definition this is part of the cohomology that corresponds to
cuspidal automorphic forms in the interpretation above.  Moreover, it is
known that $H_{\cusp}^{*}$ can only be nonzero in degrees $1$ and $2$,
and that $\dim_{\CC} H_{\cusp}^{1} = \dim_{\CC} H_{\cusp}^{2}$, via an
isomorphism $H_{\cusp}^{1} \simeq H_{\cusp}^{2}$ induced by the Hodge $*$-operator.

To summarize, we have in degrees $i=1,2$ direct sum decompositions
\[
H^{i} (\Gamma ; \CC) = H_{\cusp}^{i} (\Gamma)\oplus H_{\Eis}^{i} (\Gamma), 
\]
and the dimensions of $H_{\cusp}^{i}$ agree.  

\subsection{Cuspidal classes and Eisenstein classes}\label{ss:ccec}
Based on our data, we observed the following, denoting the class number of $\OO$ by $h$.
\begin{enumerate}
\item For $\Gamma = \GL_{2} (\OO)$, we have
\begin{itemize}
\item $\dim H^{2}_{\Eis} (\GL_{2} (\OO)) = h-1$, and $\dim H^{1}_{\Eis} (\GL_{2} (\OO)) = 0$. \cite{harder.budapest}
\item $\dim H^{*}_{\cusp} (\GL_{2} (\OO ))$ can vanish for large
$|D|$.  Indeed, for the biggest example we worked with, namely $D =
-2099$, the cuspidal cohomology vanishes.
\end{itemize}
\item For $\Gamma = \SL_{2} (\OO)$, we have
\begin{itemize}
\item If $D \not = -3, -4$, then $\dim H^{2}_{\Eis} (\SL_{2} (\OO)) =
h-1$, and $\dim H^{1}_{\Eis}
(\SL_{2} (\OO)) = h$.  \cite{serre}
\item $\dim H^{*}_{\cusp} (\SL_{2} (\OO ))$ only vanishes for finitely
many discriminants $D$.  Indeed, this is a theorem of
Grunewald--Schwermer \cite{grunewald.schwermer.1}.  Vogtmann
\cite{vogtmann} showed that the $H_{\cusp}^{1} (\SL_{2} (\OO))$ is
nonzero for only 14 discriminants, the largest being $-71$.  This is
corroborated by our computations.
\end{itemize}
\end{enumerate}

Our cuspidal data agrees with Scheutzow's $\GL_{2}/\SL_{2}$
computations \cite{scheutzow} for all the discriminants he computed.
We also have agreement with Vogtmann's $\SL_{2}$ computations, except
for two discriminants: for $D=-67$, we get $2$ compared to Vogtmann's
$3$, and for $D=-88$, we get $3$ compared to Vogtmann's $2$.

Moreover, we observe that for all $0<-D\leq 2099$,
$$\dim H^{2}_{\cusp } (\GL_{2} (\OO)) =0 \qquad \text{ only if the
ideal class group is cyclic}.$$

\begin{figure}[htb]
\begin{center}
\includegraphics[scale=0.35,angle=270,origin=c]{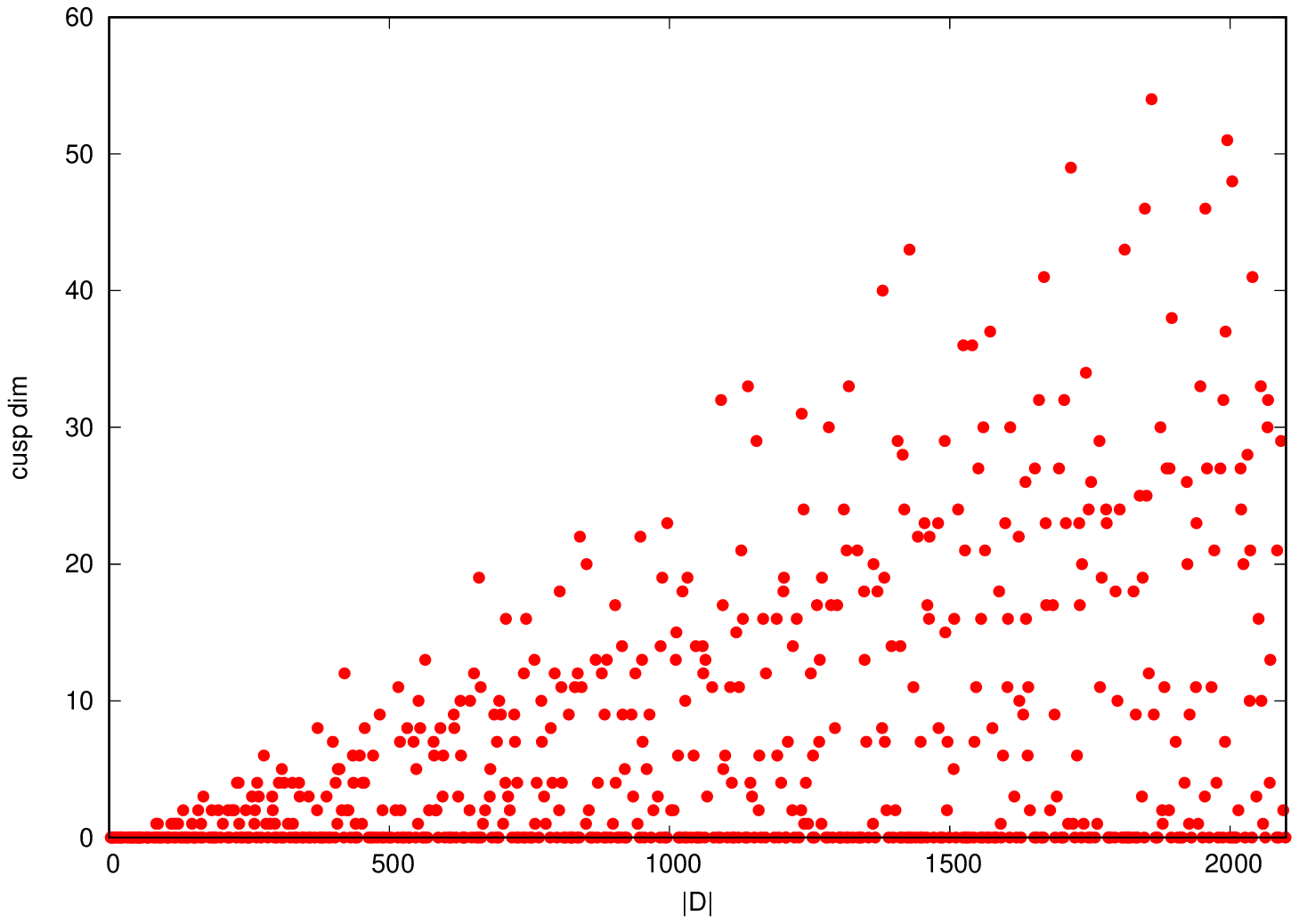}
\end{center}
\caption{$\dim H^{2}_{\cusp} (\GL_{2} (\OO_{D}))$ as a function of $|D|$.\label{fig:glcusp}}
\end{figure}

\begin{figure}[htb]
\begin{center}
\includegraphics[scale=0.35,angle=270,origin=c]{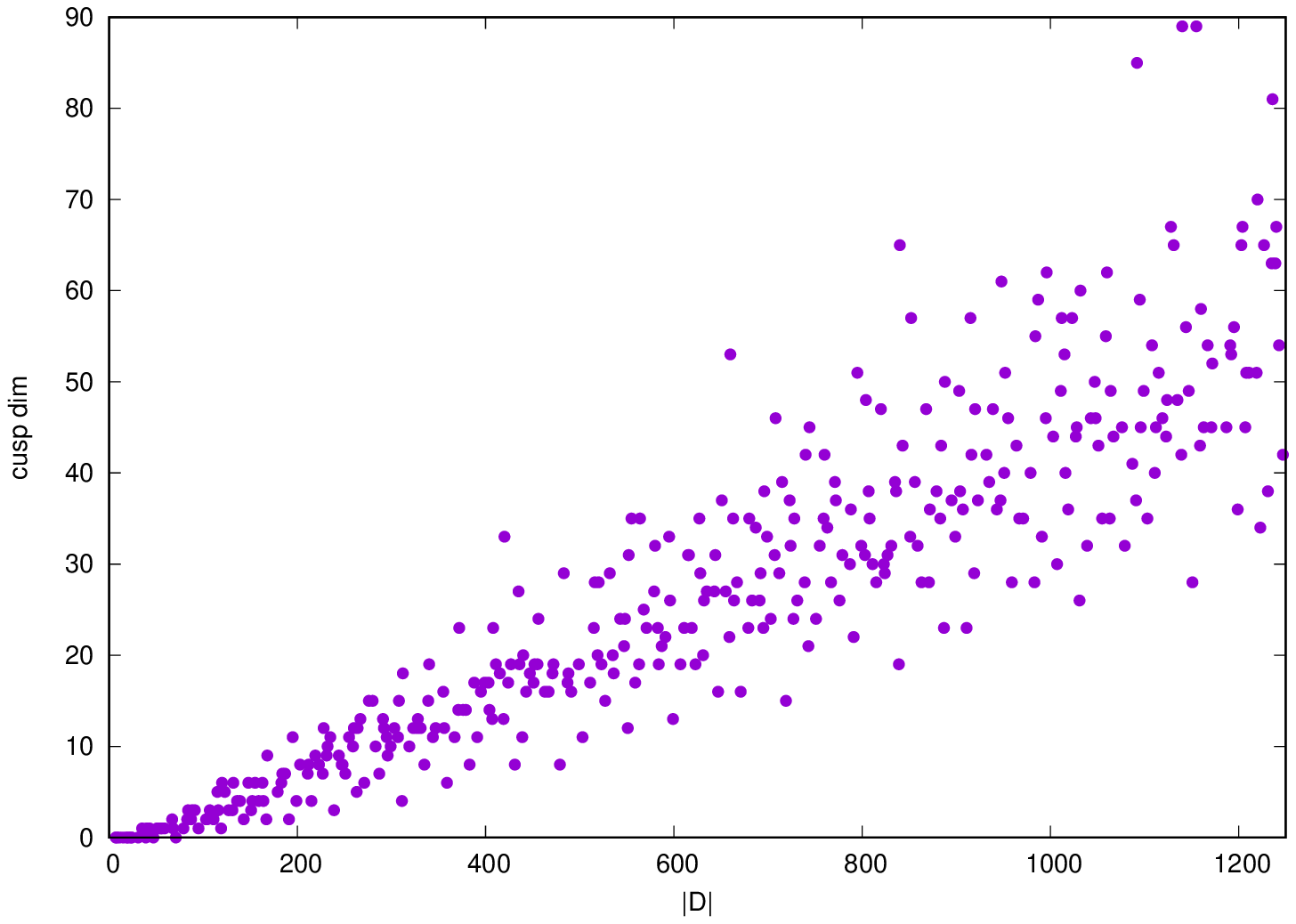}
\end{center}
\caption{$\dim H^{2}_{\cusp} (\SL_{2} (\OO_{D}))$ as a function of $|D|$.\label{fig:slcusp}}
\end{figure}

Figures~\ref{fig:glcusp} and~\ref{fig:slcusp} show $\dim H^{2}_{\cusp
}$ as a function of $|D|$. 

\subsection{Exponential growth of torsion}\label{ss:eg}
It has long been known that there can be large torsion in the homology
and cohomology of subgroups of $\GL_{2} (\OO)$ and $\SL_{2} (\OO )$.
The first examples of this phenomenon known to us (for $\PSL_{2}
(\OO_{-4})$) appear in Elstrodt--Grunewald--Mennicke \cite{egm.metz}
and Priplata's thesis \cite{priplata}.  More recently
Bergeron--Venkatesh \cite{bv}, building on ideas of Bhargava
\cite{bhargava}, formulated a precise conjecture for the growth of the
torsion in the {homology} of cocompact arithmetic subgroups
$\Gamma\subset \SL_{2} (\CC )$ as $\Gamma$ goes along a tower of
subgroups of increasing index.  For the case of congruence subgroups
of $\SL_{2} (\OO)$---which are not cocompact---analogues of the
conjecture were numerically studied by {\c{S}}eng{\"u}n \cite{sengun}
and proved (for certain coefficient modules) by Pfaff \cite{pfaff1}.

Our data supports a further variation on the theme of torsion growth,
namely for torsion growth in the degree $1$ integral Voronoi homology
of the groups $\Gamma = \GL_{2} (\OO_{D}) $.  In other words, instead of fixing a
discriminant $D$ and considering towers of congruence subgroups
$\Gamma$ of increasing index, we work at full level and let
$|D|\rightarrow \infty$.  Such computations for $\SL_{2} (\OO)$ were
carried out previously by Rahm \cite{rahm.higher}; he observed
surprisingly large torsion primes in the first integral homology of
$\SL_{2} (\OO)$, but made no conjecture about the growth.  Our data
here is for the 1st Voronoi homology group, which (away from $2$ and
$3$) is isomorphic to the 2nd integral cohomology group of $\Gamma$.
However, we believe that \emph{all} the $p$-torsion in the Voronoi
homology is arithmetically interesting, in that there should be mod
$p$ Galois representations attached to $p$-torsion Hecke eigenclasses
for all $p$ (not just $p>3$).  (For a discussion of this point in the
context of $\GL_{4} (\ZZ)$, we refer to \cite{agm6}.)  Thus its growth
behavior should be predicted by similar heuristics to
\cite{bv,bhargava}.  

\begin{figure}[htb]
\begin{center}
\includegraphics[scale=0.35,angle=270,origin=c]{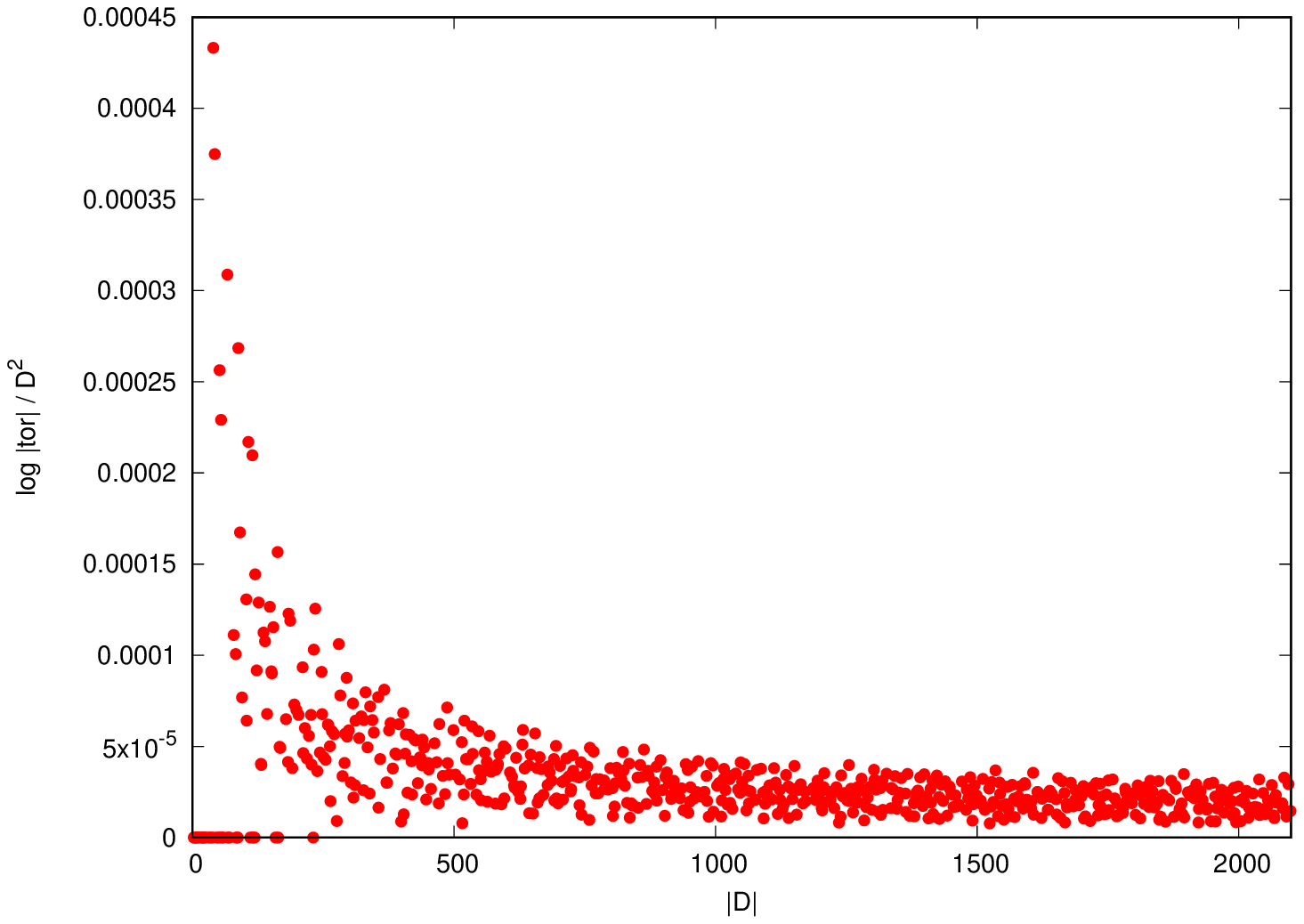}
\end{center}
\caption{The ratio $\logtor_{\GL} (D) := \log |H_{1}
(\Vor_{*} (\GL_{2} (\OO_{D})))_{\tors}|/|D|^{2}$ as a function of
$|D|$.\label{fig:gl}}
\end{figure}

\begin{figure}[htb]
\begin{center}
\includegraphics[scale=0.35,angle=270,origin=c]{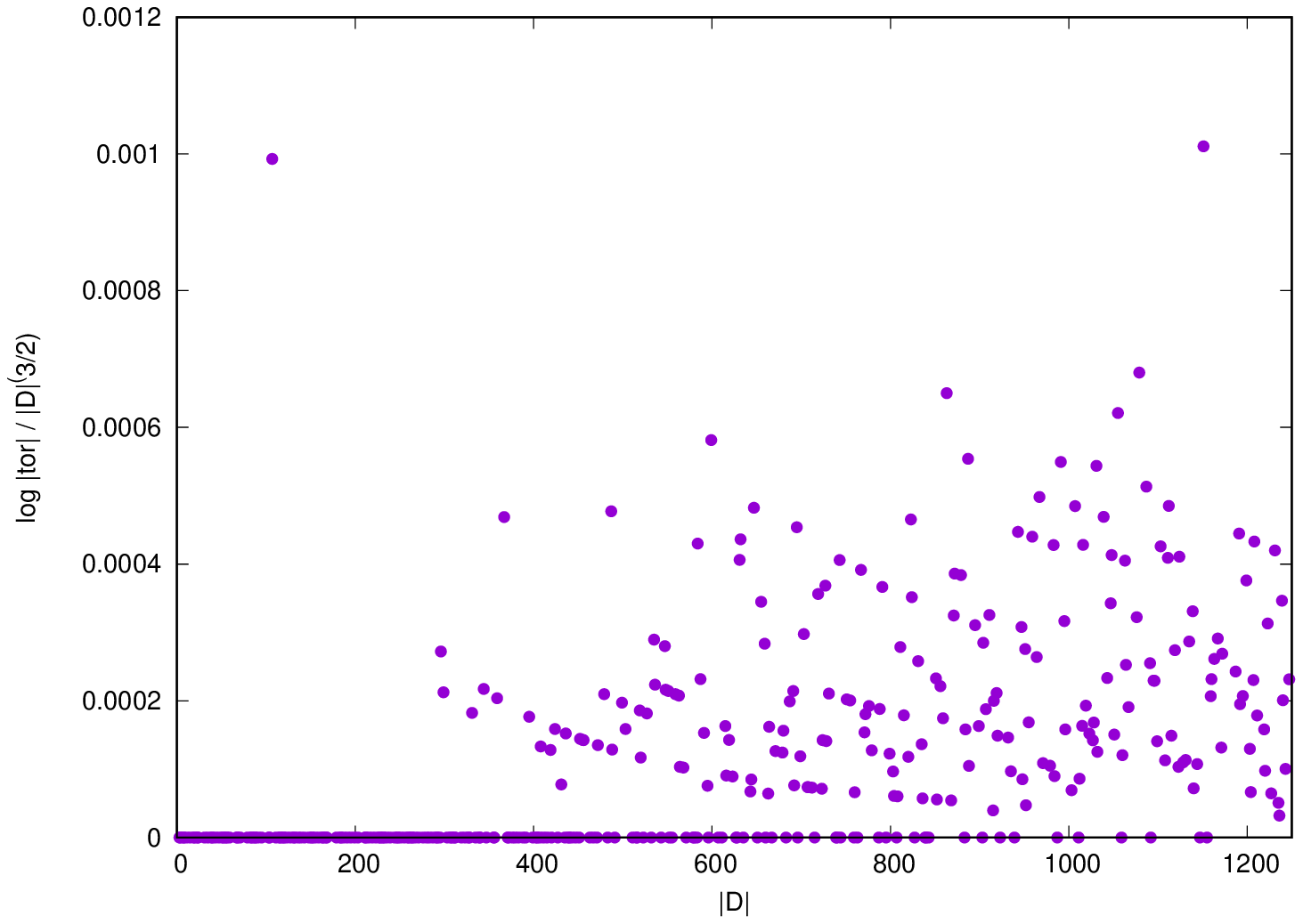}
\end{center}
\caption{The ratio $\logtor_{\SL} (D) := \log |H_{1}
(\Vor_{*} (\SL_{2} (\OO_{D})))_{\tors}|/|D|^{3/2}$ as a function of
$|D|$.\label{fig:sl}}
\end{figure}

Figure \ref{fig:gl} shows a plot of $\logtor_{\GL} (D) := \log |H_{1}
(\Vor_{*} (\GL_{2} (\OO_{D})))_{\tors}|/|D|^{2}$ as a function of
$|D|$.  It appears from this plot that the ratio $\logtor_{\GL } (D)$
tends to a nonzero limit as $|D|\rightarrow \infty$.  We note that
there is a large $2$-torsion subgroup contributing to $\logtor_{\GL }
(D)$; in fact, for the discriminants $-2035$, $-2040$ at the end of
our data range, the torsion subgroup of the Voronoi homology is
\emph{all} $2$-torsion.  Figure \ref{fig:sl} shows a similar plot for
$\logtor_{\SL} (D)$ as a function of $|D|$, although the exponent of
$|D|$ in the denominator of $\logtor$ is $3/2$, not $2$.  For this
plot, it is certainly not clear that eventually the ratio
$\logtor_{\SL} (D)$ tends to a nonzero limit.

\subsection{Rohlfs's bound and linear growth of cyclic factors}\label{ss:rbcf}

\begin{figure}[htb]
\begin{center}
\includegraphics[scale=0.35,angle=270,origin=c]{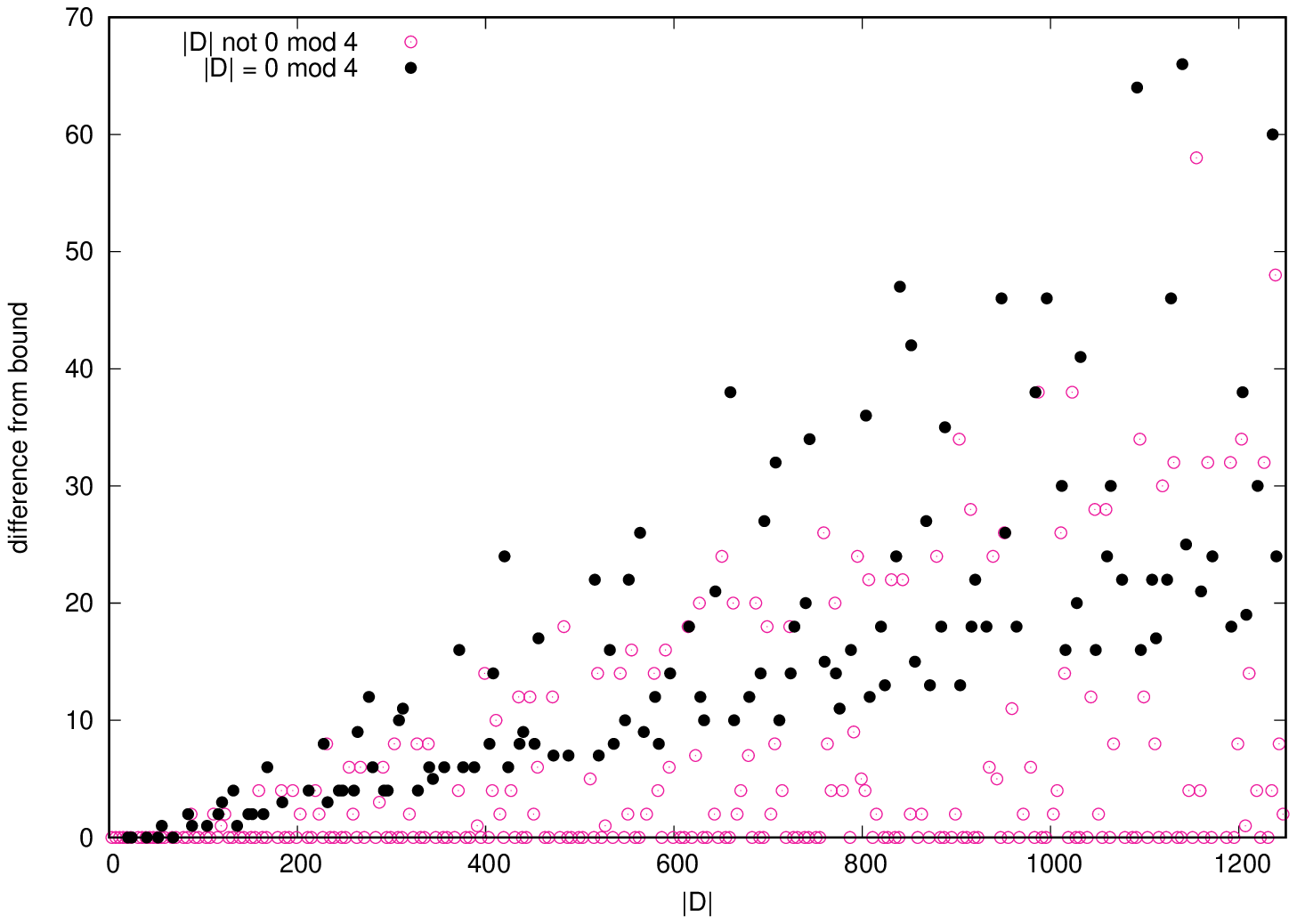}
\end{center}
\caption{The difference between the computed dimension of
$H^{1}_{\cusp }$ and Rohlfs's bound for $\SL_{2} (\OO_{D})$\label{fig:rohlfs}}
\end{figure}

In \cite{rohlfs.bound} Rohlfs gives a lower bound on $\dim
H^1_{\cusp}(\Gamma)$ where $\Gamma = \PSL(2,\OO)$.  Following an
idea of Harder in \cite{harder.budapest}, he uses that complex
conjugation $\tau$ acts on the underlying symmetric space, hence also
on the associated cohomology, giving rise to a Lefschetz fixed point
formula involving the alternating traces on the groups
$H^j$. Poincar\'e duality combined with the orientation-reversing
property of $\tau$ allows him to isolate the contribution to the trace
of $H^1_{\cusp}$, and he subsequently bounds it by analyzing the fixed
point sets of the action in detail. The main contribution of the final
bound involves Euler's totient function, which grows linearly.  It
turns out that if we replace $\PSL_2$ by $\SL_2$ and use a slightly
better bound (again due to Rohlfs, cf.~\S4.1 in loc.cit.), we find
that for $|D| \not \equiv 0 \bmod 4$ the bound is sharp in the
majority (54\%) of the computed cases.  On the other hand the gap
between the actual value and the bound for $|D| \equiv 0 \bmod 4$
appears to widen linearly with the discriminant.  We show this
behavior in Figure~\ref{fig:rohlfs}; the solid (resp.~hollow points)
correspond to $|D| \equiv$ (resp.~$\not \equiv$) $0 \bmod 4$.

\begin{figure}[htb]
\begin{center}
\includegraphics[scale=0.35,angle=270,origin=c]{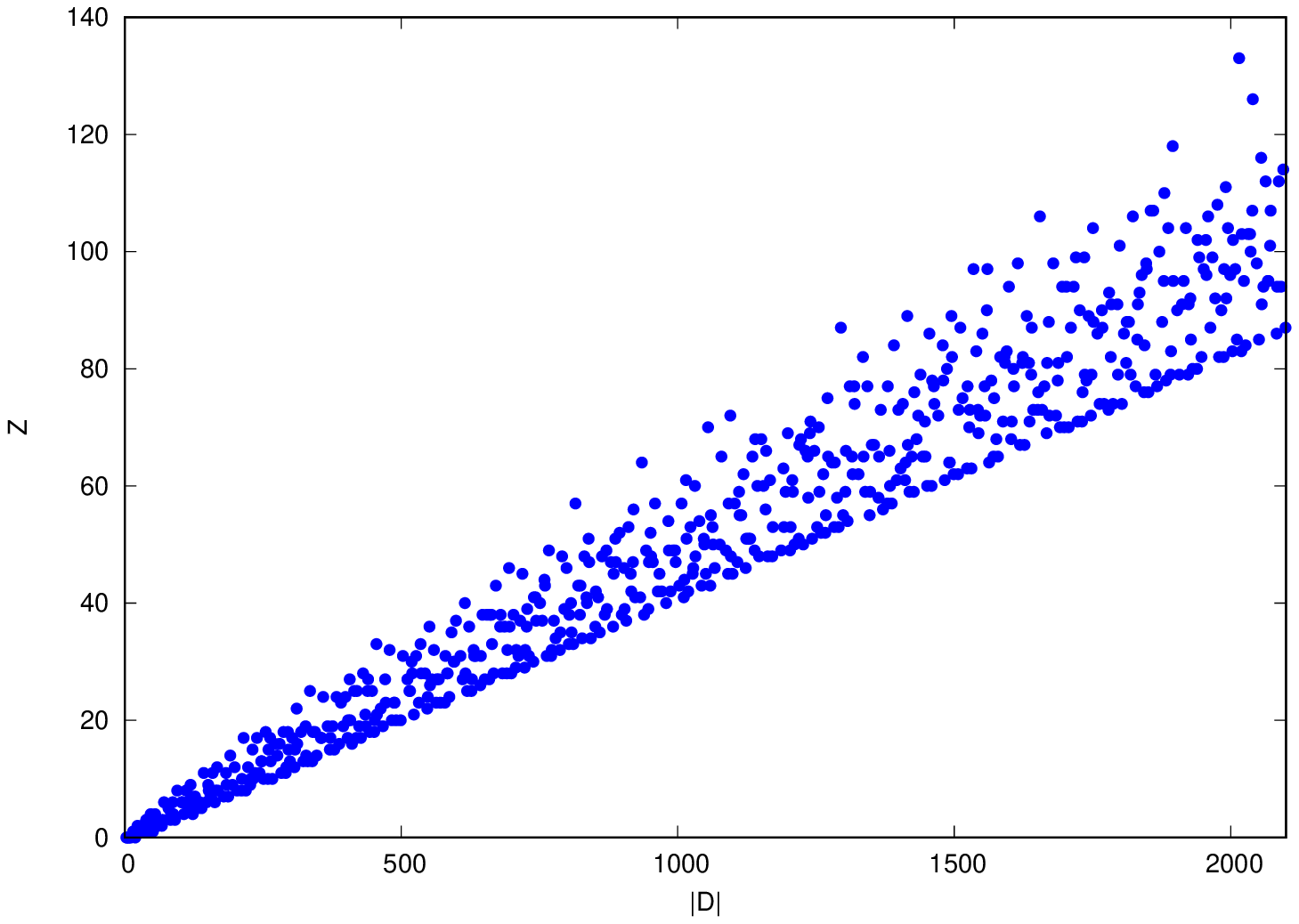}
\end{center}
\caption{The quantity $Z = Z_{D}$ as a function of
$|D|$.\label{fig:Z}}
\end{figure}

A surprisingly smooth slope arises when we pass to $\GL_2(\OO)$ and
count not only the rank of the cuspidal cohomology but rather the
number of non-trivial cyclic factors, both finite and infinite.  For
finite groups, this is often called the \emph{generator rank}.  In
particular let $Z = Z_{D}$ be the number of nontrivial cyclic factors,
both infinite and finite, of the first Voronoi homology.  This
quantity can be re-expressed as the sum of $\dim H^{2} (\GL_{2} (\OO
)) = \dim H_{\cusp}^{2} (\GL_{2} (\OO )) + \dim H_{\Eis}^{2} (\GL_{2}
(\OO ))$ plus the number of cyclic factors appearing in
$H_{1}(\Vor_{*} (\GL_{2} (\OO)))_{\tors}$.  Recall that $\dim
H_{\Eis}^{2} (\GL_{2} (\OO ))$ is one less than the class number.
Figure~\ref{fig:Z} shows a plot of this quantity as a function of
$|D|$.  One sees from this plot that the generator rank $Z$ appears to
be bounded below by a linear function of $|D|$, and that in fact
appears to meet this lower envelope infinitely often.  We have no
explanation for this phenomenon.

\section{Data}\label{sec:data}

\subsection{} The remainder of the text is devoted to presenting
details about the ranks and $p$-ranks of the homology of the Voronoi
complexes in a concise form.  We give data for each fundamental
discriminant in two tables, the first for $\GL_{2}$ and the second for
$\SL_{2}$.  The notation is as follows:

\begin{itemize}
\item The first (boldface) entry in each table cell is the
discriminant of the imaginary quadratic field $F=\QQ(\sqrt{D})$.
\item The second entry gives the cyclic factors of the class group of
$F$.
\item Next we give two entries encoding the Voronoi homology group
$H_{1}$; modulo the Serre class $\cS_{3}$ this is the same as the
integral group cohomology $H^{2}$.  The third entry is the rank of the
cuspidal subspace; the full rank of $H^{2}$ can be recovered from this
by adding $h-1$, where $h$ is the class number of ${F}$.
\item Finally the fourth entry gives the invariant factors of the
torsion subgroup of the Voronoi homology $H_{1}$.  Since there are
many invariant factors, we have used exponential notation: the
notation $(d_{1}^{n_{1}},\dotsc ,d_{k}^{n_{k}})$ means that the
torsion subgroup is
\[
\prod_{i=1}^{k}(\ZZ/d_{i}\ZZ )^{n_{i}}. 
\]
Away from the $2$- and $3$-torsion, this is the same as the torsion
subgroup of $H^{2}$.
\end{itemize}
We omit the Voronoi homology group $H_{2}$, which modulo the Serre
class is the same as the integral group cohomology $H^{1}$, since it
is easily recoverable from the data for $H^{2}$. In particular the
cuspidal rank agrees with that for $H^{2}$, and for $\GL_{2}$ this is
full rank of $H^{1}$.  For $\SL_{2}$ and $D<-4$ the full rank is the
cuspidal rank plus $h$ (for $D=-3, -4$ the rank of $H^{1}$ vanishes).
We omit the invariant factors for $H_{1}$ since in all cases they were annihilated
by $12$.

\subsection{} As an example, consider the $\GL_{2}$ table entry for
discriminant $-1007$.  The notation ${(30)}$ means the
class group of $\OO_{-1007}$ is $\ZZ /30\ZZ$, and thus $h=30$.  The
next two entries are 
\[
2\quad \quad (2^{24}, 534, 1602).
\]
This means $\dim H^{2}_{\cusp} (\GL_{2} (\OO_{-1007}); \CC)$
is $2$.  Since $\dim H^{2}_{\Eis} (\GL_{2} (\OO_{-1007} ); \CC) = h-1 = 29$,
the full Voronoi homology group $H_{1} (\GL_{2} (\OO_{-1007}))$ is
isomorphic to 
\[
\ZZ^{31} \times (\ZZ /2\ZZ)^{24}\times \ZZ /534\ZZ \times \ZZ /1602\ZZ
\simeq \ZZ^{31} \times (\ZZ /2\ZZ)^{26}\times (\ZZ/3\ZZ )^{3}\times
(\ZZ / 89\ZZ)^{2}. 
\]
We note that the torsion contains the
large prime $89$.

\bibliographystyle{amsplain_initials_eprint}
\bibliography{gl2imag}

\newpage

\iflandscapetable
\begin{landscape}
\fi

\thispagestyle{empty}
\pagestyle{empty}

\section{\texorpdfstring{Cohomology data for $\GL_{2} (\OO) $}{Cohomology data for $\GL_{2}(\OO) $}} \label{sec:gl2data}
\vspace{-0.35cm}

\begin{center}
\footnotesize

\iflandscapetable

Cohomology data for $\GL_2(\OO)$.
\end{center}
\newpage

\begin{center}
\footnotesize

\iflandscapetable
%
Cohomology data for $\GL_2(\OO)$.
\end{center}
\newpage

\begin{center}
\footnotesize

\iflandscapetable
%
Cohomology data for $\GL_2(\OO)$.
\end{center}
\newpage

\begin{center}
\footnotesize

\iflandscapetable
%
Cohomology data for $\GL_2(\OO)$.
\end{center}
\newpage

\begin{center}
\footnotesize

\iflandscapetable
%
Cohomology data for $\GL_2(\OO)$.
\end{center}
\newpage

\begin{center}
\footnotesize

\iflandscapetable
%
Cohomology data for $\GL_2(\OO)$.
\end{center}
\newpage

\begin{center}
\footnotesize

\iflandscapetable
%
Cohomology data for $\GL_2(\OO)$.
\end{center}
\newpage

\newpage

\section{\texorpdfstring{Cohomology data for $\SL_{2} (\OO) $}{Cohomology data for $\SL_{2}(\OO) $}} \label{sec:sl2data}
\vspace{-0.35cm}

\begin{center}
\footnotesize

\iflandscapetable
%
Cohomology data for $\SL_2(\OO)$.
\end{center}
\newpage

\begin{center}
\footnotesize

\iflandscapetable
%
Cohomology data for $\SL_2(\OO)$.
\end{center}
\newpage

\begin{center}
\footnotesize

\iflandscapetable
%
Cohomology data for $\SL_2(\OO)$.
\end{center}
\newpage

\begin{center}
\footnotesize

\iflandscapetable
\begin{longtable}{cccc|cccc|cccc}
\fi

\toprule
\midrule
\endfirsthead
\toprule
\midrule
\endhead
\bottomrule
\endfoot
\bottomrule
\endlastfoot
$D$ & $\Cl(D)$ & $H^2_{\cusp}$ & $(H_1)_{\tors}$ &$D$ & $\Cl(D)$ & $H^2_{\cusp}$ & $(H_1)_{\tors}$ &$D$ & $\Cl(D)$ & $H^2_{\cusp}$ & $(H_1)_{\tors}$ \\
\hline
$\mathbf{-1199}$ & ${(38)}$ & $36$ & $(2461^{2})$ &$\mathbf{-1203}$ & ${(6)}$ & $65$ & $(15^{2})$ &$\mathbf{-1204}$ & ${(2,4)}$ & $67$ & $(2^{4})$ \\
\hline
$\mathbf{-1207}$ & ${(18)}$ & $45$ & $(125^{2})$ &$\mathbf{-1208}$ & ${(12)}$ & $51$ & $(2^{14},40,120)$ &$\mathbf{-1211}$ & ${(14)}$ & $51$ & $(43^{2})$ \\
\hline
$\mathbf{-1219}$ & ${(6)}$ & $51$ & $(29^{2})$ &$\mathbf{-1220}$ & ${(2,8)}$ & $70$ & $(2^{6})$ &$\mathbf{-1223}$ & ${(35)}$ & $34$ & $(808^{2})$ \\
\hline
$\mathbf{-1227}$ & ${(4)}$ & $65$ & $(4^{2})$ &$\mathbf{-1231}$ & ${(27)}$ & $38$ & $(2^{8},544^{2})$ &$\mathbf{-1235}$ & ${(2,6)}$ & $63$ & $(3^{2})$ \\
\hline
$\mathbf{-1236}$ & ${(2,6)}$ & $81$ & $(2^{2})$ &$\mathbf{-1239}$ & ${(2,16)}$ & $63$ & $(7^{2},273^{2})$ &$\mathbf{-1240}$ & ${(2,4)}$ & $67$ & $(2^{6},10^{2})$ \\
\hline
$\mathbf{-1243}$ & ${(4)}$ & $54$ & $(3^{4})$ &$\mathbf{-1247}$ & ${(26)}$ & $42$ & $(163^{2})$ \\
\hline
\end{longtable}
Cohomology data for $\SL_2(\OO)$.
\end{center}
\newpage

\iflandscapetable
\end{landscape}
\fi

\newpage

\end{document}